\documentclass[12pt]{article}
\usepackage{latexsym}
\usepackage{amsfonts}
\usepackage{amsmath}
\usepackage{amssymb}
\usepackage{graphicx}
\usepackage{latexsym}
\usepackage{amsfonts}
\usepackage{graphicx}
\usepackage{psfrag}

\input{tcilatex}
\begin{document}

\qquad 

\thispagestyle{empty}

\begin{center}
{\Large \textbf{\ Complementary upper bounds for fourth central moment with
extensions and applications }}

\vskip0.5inR. Sharma, R. Kumar, R. Saini and G. Kapoor\\[0pt]
Department of Mathematics\\[0pt]
H. P. University, Shimla -5,\\[0pt]
India - 171 005\\[0pt]
email: rajesh.sharma.hpn@nic.in
\end{center}

\bigskip \vskip1.5in \textbf{Abstract. }We prove some inequalities involving
fourth central moment of a random variable that takes values in a given
finite interval. Both discrete and continuous cases are considered. Bounds
for the spread are obtained when a given $n\times n$ complex matrix has real
eigenvalues. Likewise, we discuss bounds for the spans of polynomial
equations.

\vskip0.5in \textbf{AMS classification \quad } 60E15, 15A42, 12D10


\bigskip

\bigskip

\bigskip

\bigskip

\bigskip

\bigskip

\bigskip

\bigskip

\section{Introduction}

\setcounter{equation}{0}The $r^{th}\ $central moment $\mu _{r}$ of a random
variable $X$ in $\left[ m,M\right] $ for the continuous and discrete cases
respectively are defined as%
\begin{equation}
\mu _{r}=\dint\limits_{m}^{M}\left( x-\mu _{1}^{^{\prime }}\right)
^{r}f\left( x\right) dx\text{ \ or \ }\mu _{r}=\sum_{i=1}^{n}p_{i}\left(
x_{i}-\mu _{1}^{^{\prime }}\right) ^{r},  \tag{1.1}
\end{equation}%
where%
\begin{equation}
\mu _{1}^{^{\prime }}=\dint\limits_{m}^{M}xf\left( x\right) \text{ or \ }\mu
_{1}^{^{\prime }}=\sum_{i=1}^{n}p_{i}x_{i},  \tag{1.2}
\end{equation}%
\ $f\left( x\right) $ and $p_{i}$ are corresponding probability densities
and probability functions such that%
\begin{equation}
\dint\limits_{m}^{M}f\left( x\right) =1\text{ \ or}\sum_{i=1}^{n}p_{i}=1. 
\tag{1.3}
\end{equation}%
We denote by $m_{r}$ the $r^{th}\ $central moment of $n$ real numbers $%
x_{1},x_{2},...,x_{n},$%
\begin{equation}
m_{r}=\frac{1}{n}\sum_{i=1}^{n}\left( x_{i}-m_{1}^{^{\prime }}\right) ^{r}, 
\tag{1.4}
\end{equation}%
where $m_{1}^{^{\prime }}=\frac{1}{n}\underset{i=1}{\overset{n}{\sum }}x_{i}$
is the arithmetic mean.

Bounds on the variance ($\sigma ^{2}=\mu _{2},\ S^{2}=m_{2}$), their
extensions and applications have been studied extensively in literature; see 
\cite[5, 13-19]{3}. The well-known Popoviciu inequality gives an upper bound
for the variance of a random variable, \cite{13}, 
\begin{equation}
\mu _{2}\leq \frac{\left( M-m\right) ^{2}}{4}.  \tag{1.5}
\end{equation}%
Nagy's inequality \cite{11} provides a complementary lower bound for the
variance of $n$ real numbers $x_{i};$ $m\leq x_{i}\leq M,\ i=1,2,...n,$ 
\begin{equation}
S^{2}\geq \frac{\left( M-m\right) ^{2}}{2n}.  \tag{1.6}
\end{equation}%
Such inequalities are also useful in many other contexts. For example,
Wolkowicz and Styan \cite{19} have observed that if the eigenvalues$\ $of an 
$n\times n$ complex matrix are all real, as in case of Hermitian matrices,
the inequalities (1.5) and (1.6) provide bounds for the spread of a matrix,
spd$\left( A\right) =\underset{i,j}{\max }\left\vert \lambda _{i}-\lambda
_{j}\right\vert .\ $Let $B=A-\frac{\text{tr}A}{n}I,\ $where tr$A$ denotes
the trace of $A$. Then,%
\begin{equation}
\frac{4}{n}\text{tr}B^{2}\leq \text{spd}\left( A\right) ^{2}\leq 2\text{tr}%
B^{2}.  \tag{1.7}
\end{equation}

Further, let $\mathbb{M}(n)$ denotes the C$^{\ast }-$algebra of all $n\times
n$ complex matrices and let $\Phi :\mathbb{M}(n)\rightarrow \mathbb{M}(k)$
be a positive unital linear map \cite{4}. The inequality of Bhatia and Davis 
\cite{3} says that if the spectrum of a Hermitian matrix $A$ is contained in
the interval $\left[ m,M\right] $, then%
\begin{equation}
\Phi \left( A^{2}\right) -\Phi \left( A\right) ^{2}\leq \frac{\left(
M-m\right) ^{2}}{4}=\frac{\text{spd}\left( A\right) ^{2}}{4},  \tag{1.8}
\end{equation}%
for every positive unital linear map $\Phi .$ This gives a noncommutative
analogue of the inequality (1.5) and yields many old and new bounds for the
spread of a matrix. This is demonstrated in \cite{5}.

Likewise, the inequalities (1.5) and (1.6) provide bounds for the span of
polynomial, see \cite{14} and Section 4, below.

Such basic inequalities, their further refinements, extensions and
alternative proofs have been studied by several authors. In particular,
Sharma et al \cite[16]{15} have proved that%
\begin{equation}
\frac{\mu _{2}^{2}-\left( \mu _{1}^{^{\prime }}-m\right) ^{2}\mu _{2}}{\mu
_{1}^{^{\prime }}-m}\leq \mu _{3}\leq \frac{\left( M-\mu _{1}^{^{\prime
}}\right) ^{2}\mu _{2}-\mu _{2}^{2}}{M-\mu _{1}^{^{\prime }}}  \tag{1.9}
\end{equation}%
and%
\begin{equation}
\sigma ^{2}+\left( \frac{\mu _{3}}{2\sigma ^{2}}\right) ^{2}\leq \frac{%
\left( M-m\right) ^{2}}{4}.  \tag{1.10}
\end{equation}%
The inequality (1.10) provides a refinement of the Popoviciu inequality
(1.5). The inequalities (1.9) and (1.10) yield bounds for the eigenvalues
and spread of a Hermitian matrix. Likewise, these inequalities provide
bounds for the roots of polynomial equations. See \cite{15}.

We focus here on inequalities involving fourth central moment ($\mu _{4}$).
One such inequality in literature is Pearson's inequality \cite{12} which
gives an interesting relation between two important parameters of
statistical distributions namely skewness ($\alpha _{3}$) and kurtosis ($%
\alpha _{4}$),%
\begin{equation*}
\alpha _{4}\geq 1+\alpha _{3}^{2}\ ,
\end{equation*}%
where 
\begin{equation}
\alpha _{3}=\sqrt{\frac{m_{3}^{2}}{m_{2}^{3}}}\ \ \text{and\ \ }\alpha _{4}=%
\frac{m_{4}}{m_{2}^{2}}.  \tag{1.11}
\end{equation}%
For more details, see \cite[17]{16} and references therein.

We derive some inequalities involving fourth central moment and discuss
related extensions and applications. We prove an analogue of the Popoviciu
inequality (1.5) for fourth central moment (Theorem 2.1, below). Our main
result (Theorem 2.2) gives bounds for the fourth central moment in terms of
second and third central moments. The inequalities involving first four
central moments and range of the random variable are obtained (Corollary
2.3-2.4). This also provides a relation among skewness, kurtosis and
studentized range (Corollary 2.5). It is shown that the inequality (1.10)
provides a refinement of the inequality for third central moment in terms of
the range of the random variable (Theorem 2.3). A generalization of the Nagy
inequality (1.6) is proved for the $s^{th}$ central moment, $s=2r$ (Theorem
2.4). We obtain bounds for the spread of a Hermitian matrix (Theorem 3.1 and
3.3). Likewise, bounds for the span of polynomial are discussed (Theorem
4.1-4.2).

\section{Main results}

\setcounter{equation}{0}It is enough to prove the following results for the
case when $X$ is a discrete random variable taking finitely many values $%
x_{1},x_{2},...,x_{n}$ with probabilities $p_{1},p_{2},...,p_{n}\ ,\ $%
respectively. The arguments are similar for the case when $X$ is a
continuous random variable.

\vskip0.1in\textbf{Theorem 2.1. }Let $X$ be a discrete or continuous random
variable taking values in\textbf{\ }$\left[ m,M\right] .$ Then%
\begin{equation}
\mu _{4}\leq \frac{\left( M-m\right) ^{4}}{12}.  \tag{2.1}
\end{equation}

\vskip0.1in\textbf{Proof. }For $\alpha \leq y\leq \beta ,$ we have%
\begin{equation}
\left( y-\alpha \right) \left( y-\beta \right) \left( \left( y+\frac{\alpha
+\beta }{2}\right) ^{2}+\frac{\alpha ^{2}+\beta ^{2}+\left( \alpha +\beta
\right) ^{2}}{4}\right) \leq 0.  \tag{2.2}
\end{equation}%
Put $y=x_{i}-\mu _{1}^{^{\prime }},\ \alpha =m-\mu _{1}^{^{\prime }}\ $and $%
\beta =M-\mu _{1}^{^{\prime }}$ in (2.2), multiply both sides by $p_{i}$,
add $n$ inequalities, $i=1,2,...n$ and use (1.1)-(1.3), we see that%
\begin{equation}
\mu _{4}\leq \left( \mu _{1}^{^{\prime }}-m\right) \left( M-\mu
_{1}^{^{\prime }}\right) \left( \left( \mu _{1}^{^{\prime }}-m\right)
^{2}+\left( M-\mu _{1}^{^{\prime }}\right) ^{2}-\left( \mu _{1}^{^{\prime
}}-m\right) \left( M-\mu _{1}^{^{\prime }}\right) \right) .  \tag{2.3}
\end{equation}%
The inequality (2.1) now follows from (2.3) and the fact that the function%
\begin{equation*}
h(x)=\left( x-m\right) \left( M-x\right) \left( \left( x-m\right)
^{2}+\left( M-x\right) ^{2}-\left( x-m\right) \left( M-x\right) \right) ,
\end{equation*}%
achieves its maximum at%
\begin{equation*}
x=\frac{m+M}{2}\pm \frac{m-M}{2\sqrt{3}},
\end{equation*}%
where%
\begin{equation*}
h(x)\leq \frac{\left( M-m\right) ^{4}}{12}.\ \ \blacksquare \ 
\end{equation*}%
The sign of equality holds in (2.3) if and only if $n=2$. In this case, $%
m_{4}=\frac{\left( M-m\right) ^{4}}{16}$. Equality holds in (2.1) for $n=2;\
x_{1}=m\ $and $x_{2}=M$ with$\ p_{1}=\frac{1}{2}\pm \frac{1}{2\sqrt{3}}\ $%
and $p_{2}=\frac{1}{2}\mp \frac{1}{2\sqrt{3}}$.

Pearson's inequality \cite{12} gives a lower bound for the fourth central
moment,%
\begin{equation}
\mu _{4}\geq \frac{\mu _{3}^{2}}{\mu _{2}}+\mu _{2}^{2}.  \tag{2.4}
\end{equation}%
We derive a complementary upper bound in the following theorem.

\vskip0.1in\textbf{Theorem 2.2. }Let $X$ be a discrete or continuous random
variable taking values in\textbf{\ }$\left[ m,M\right] .\ $Then, 
\begin{equation}
\mu _{4}\leq \left( \mu _{1}^{^{\prime }}-m\right) \left( M-\mu
_{1}^{^{\prime }}\right) \mu _{2}+\left( m+M-2\mu _{1}^{^{\prime }}\right)
\mu _{3}-\frac{\left( \mu _{3}-\left( m+M-2\mu _{1}^{^{\prime }}\right) \mu
_{2}\right) ^{2}}{\left( \mu _{1}^{^{\prime }}-m\right) \left( M-\mu
_{1}^{^{\prime }}\right) -\mu _{2}},  \tag{2.5}
\end{equation}%
where$\ \mu _{2}\neq \left( \mu _{1}^{^{\prime }}-m\right) \left( M-\mu
_{1}^{^{\prime }}\right) .$

\vskip0.1in\textbf{Proof. }Let $\alpha \leq y\leq \beta .$ Then, for any
real number $\gamma ,$%
\begin{equation}
\left( y-\alpha \right) \left( y-\beta \right) \left( y-\gamma \right)
^{2}\leq 0.  \tag{2.6}
\end{equation}%
Put $y=x_{i}-\mu _{1}^{^{\prime }},\ \alpha =m-\mu _{1}^{^{\prime }}\ $and $%
\beta =M-\mu _{1}^{^{\prime }}$ in (2.6), multiply both sides by $p_{i}$,
add $n$ inequalities, $i=1,2,...n$ and use (1.1)-(1.3), we get%
\begin{eqnarray}
\mu _{4} &\leq &\left( \left( \mu _{1}^{^{\prime }}-m\right) \left( M-\mu
_{1}^{^{\prime }}\right) -\mu _{2}\right) \gamma ^{2}-2\left( \left(
m+M-2\mu _{1}^{^{\prime }}\right) \mu _{2}-\mu _{3}\right) \gamma  \notag \\
&&+\left( m+M-2\mu _{1}^{^{\prime }}\right) \mu _{3}+\left( \mu
_{1}^{^{\prime }}-m\right) \left( M-\mu _{1}^{^{\prime }}\right) \mu _{2}\ \
.  \TCItag{2.7}
\end{eqnarray}%
The inequality (2.7) is valid for every real number $\gamma $ and gives
least upper bound for%
\begin{equation}
\gamma =\frac{\left( m+M-2\mu _{1}^{^{\prime }}\right) \mu _{2}-\mu _{3}}{%
\left( \mu _{1}^{^{\prime }}-m\right) \left( M-\mu _{1}^{^{\prime }}\right)
-\mu _{2}}.  \tag{2.8}
\end{equation}%
Substitute the value of $\gamma $\ from (2.8) in (2.7); a little calculation
leads to (2.5).\ \ \ $\blacksquare $

Note that $\mu _{2}=\left( \mu _{1}^{^{\prime }}-m\right) \left( M-\mu
_{1}^{^{\prime }}\right) $ if and only if every $x_{i}$ is equal either to $%
m $ or to $M$, see [3]. So, (2.5) is not valid for $n=2$. Equality holds in
(2.5) when%
\begin{equation*}
x_{1}=m,\ x_{2}=\frac{\left( m+M-2\mu _{1}^{^{\prime }}\right) \mu _{2}-\mu
_{3}}{\left( \mu _{1}^{^{\prime }}-m\right) \left( M-\mu _{1}^{^{\prime
}}\right) -\mu _{2}}\text{ \ and }x_{3}=M;\ n=3.
\end{equation*}%
Pearson's inequality (2.4) implies that $\mu _{2}\mu _{4}-\mu _{2}^{3}-\mu
_{3}^{2}\geq 0.\ $We prove a complementary upper bound in the following
theorem.

\vskip0.1in\textbf{Corollary 2.3. }Under the conditions of Theorem 2.2, we
have%
\begin{equation}
\mu _{2}\mu _{4}-\mu _{2}^{3}-\mu _{3}^{2}\leq \frac{\left( \left( \mu
_{1}^{^{\prime }}-m\right) \left( M-\mu _{1}^{^{\prime }}\right) \left(
M-m\right) \right) ^{2}}{27}\leq \frac{\left( M-m\right) ^{6}}{432}. 
\tag{2.9}
\end{equation}

\vskip0.1in\textbf{Proof. }From (2.5 ), we have%
\begin{eqnarray}
\mu _{2}\mu _{4}-\mu _{2}^{3}-\mu _{3}^{2} &\leq &\left( \mu _{1}^{^{\prime
}}-m\right) \left( M-\mu _{1}^{^{\prime }}\right) \mu _{2}^{2}+\left(
m+M-2\mu _{1}^{^{\prime }}\right) \mu _{2}\mu _{3}  \notag \\
&&+\frac{\mu _{2}\left( \mu _{3}-\left( m+M-2\mu _{1}^{^{\prime }}\right)
\mu _{2}\right) ^{2}}{\mu _{2}-\left( \mu _{1}^{^{\prime }}-m\right) \left(
M-\mu _{1}^{^{\prime }}\right) }-\mu _{2}^{3}-\mu _{3}^{2}\ .  \TCItag{2.10}
\end{eqnarray}%
One can easily see, on using derivatives that the right hand side expression
in (2.10) is maximum at%
\begin{equation*}
\mu _{3}=\frac{1}{2}\frac{\left( m+M-2\mu _{1}^{^{\prime }}\right) \mu _{2}}{%
\left( \mu _{1}^{^{\prime }}-m\right) \left( M-\mu _{1}^{^{\prime }}\right) }%
\left( \mu _{2}+\left( \mu _{1}^{^{\prime }}-m\right) \left( M-\mu
_{1}^{^{\prime }}\right) \right) .
\end{equation*}%
So,%
\begin{equation}
\mu _{2}\mu _{4}-\mu _{2}^{3}-\mu _{3}^{2}\leq \left( 1-\frac{\mu _{2}}{%
\left( \mu _{1}^{^{\prime }}-m\right) \left( M-\mu _{1}^{^{\prime }}\right) }%
\right) \left( \frac{M-m}{2}\right) ^{2}\mu _{2}^{2}.  \tag{2.11}
\end{equation}%
The first inequality (2.9) now follows from (2.11) and the fact that right
hand side expression (2.11) is maximum at%
\begin{equation*}
\mu _{2}=\frac{2}{3}\left( \mu _{1}^{^{\prime }}-m\right) \left( M-\mu
_{1}^{^{\prime }}\right) .
\end{equation*}%
Using arithmetic-geometric mean inequality, we have%
\begin{equation}
\left( \mu _{1}^{^{\prime }}-m\right) \left( M-\mu _{1}^{^{\prime }}\right)
\leq \left( \frac{M-m}{2}\right) ^{2}.  \tag{2.12}
\end{equation}%
The second inequality (2.9) follows from (2.12). \ \ $\blacksquare $

We now prove one more inequality complementary to Pearson's inequality $\mu
_{4}-\mu _{2}^{2}-\frac{\mu _{3}^{2}}{\mu _{2}}\geq 0$ in the following
theorem.

\vskip0.1in\textbf{Corollary 2.4. }Under the conditions of Theorem 2.2, we
have%
\begin{equation}
\mu _{4}-\mu _{2}^{2}-\frac{\mu _{3}^{2}}{\mu _{2}}\leq \left( \mu
_{1}^{^{\prime }}-m\right) \left( M-\mu _{1}^{^{\prime }}\right) \left( 
\frac{M-m}{4}\right) ^{2}\leq \frac{\left( M-m\right) ^{4}}{64}.  \tag{2.13}
\end{equation}%
\vskip0.1in\textbf{Proof. }From (2.11), we have%
\begin{equation}
\mu _{4}-\mu _{2}^{2}-\frac{\mu _{3}^{2}}{\mu _{2}}\leq \left( 1-\frac{\mu
_{2}}{\left( \mu _{1}^{^{\prime }}-m\right) \left( M-\mu _{1}^{^{\prime
}}\right) }\right) \left( \frac{M-m}{2}\right) ^{2}\mu _{2}.  \tag{2.14}
\end{equation}%
The first inequality (2.13) follows from the fact the right hand side
expression in (2.14) is maximum at%
\begin{equation*}
\mu _{2}=\frac{1}{2}\left( \mu _{1}^{^{\prime }}-m\right) \left( M-\mu
_{1}^{^{\prime }}\right) .
\end{equation*}%
The second inequality (2.13) follows from (2.12). \ \ $\blacksquare $

The studentized range $q$ of $n$ real numbers $x_{i};$ $m\leq x_{i}\leq M,\
i=1,2,...n$ is defined as%
\begin{equation}
q=\frac{M-m}{S},  \tag{2.15}
\end{equation}%
where $S$ is standard deviation.$\ $We now find an interesting relation
among studentized range, skewness and kurtosis.

\vskip0.1in\textbf{Corollary 2.5. }For $m\leq x_{i}\leq M,\ i=1,2,...n,$ we
have%
\begin{equation}
\alpha _{4}-\alpha _{3}^{2}\leq \frac{q^{2}}{4},  \tag{2.16}
\end{equation}%
where $\alpha _{3},\alpha _{4}$ and $q$ are respectively defined by (1.11)
and (2.15).

\vskip0.1in\textbf{Proof. }Divide both sides of (2.11) by $\mu _{2}^{3}$, we
see that%
\begin{equation}
\frac{\mu _{4}}{\mu _{2}^{2}}-\frac{\mu _{3}^{2}}{\mu _{2}^{3}}\leq \frac{%
\left( M-m\right) ^{2}}{4\mu _{2}}+\left( 1-\frac{\left( M-m\right) ^{2}}{%
4\left( \mu _{1}^{^{\prime }}-m\right) \left( M-\mu _{1}^{^{\prime }}\right) 
}\right) .  \tag{2.17}
\end{equation}%
Combine (2.12) and (2.17), we get that%
\begin{equation}
\frac{\mu _{4}}{\mu _{2}^{2}}-\frac{\mu _{3}^{2}}{\mu _{2}^{3}}\leq \frac{%
\left( M-m\right) ^{2}}{4\mu _{2}}.  \tag{2.18}
\end{equation}%
The inequality (2.18) implies (2.16), use (1.11) and (2.15)$.\ \ \
\blacksquare $

\textbf{Remark. }The $r^{th}$\ order moment about origin is defined as%
\begin{equation*}
\mu _{r}^{^{\prime }}=\dint\limits_{m}^{M}x^{r}f\left( x\right) dx\text{ \
or \ }\mu _{r}^{^{\prime }}=\sum_{i=1}^{n}p_{i}x_{i}^{r}.
\end{equation*}%
On using the well-known relations, $\mu _{2}=\mu _{2}^{^{\prime }}-\mu
_{1}^{\prime ^{2}},\ \mu _{3}=\mu _{3}^{^{\prime }}-3\mu _{1}^{^{\prime
}}\mu _{2}^{\prime }+2\mu _{1}^{\prime ^{3}}$ and $\mu _{4}=\mu
_{4}^{^{\prime }}-4\mu _{1}^{^{\prime }}\mu _{3}^{\prime }+6\mu _{1}^{\prime
^{2}}\mu _{2}^{^{\prime }}-3\mu _{1}^{\prime ^{2}}$ in above inequalities,
we can write the inequalities involving moments about origin of discrete and
continuous distributions. For example, the inequalities (2.5), (2.9) and
(2.13) respectively give%
\begin{equation*}
\mu _{4}^{^{\prime }}\leq \left( m+M\right) \mu _{3}^{^{\prime }}-mM\mu
_{2}^{^{\prime }}-\frac{\left( \mu _{3}^{^{\prime }}-\left( m+M\right) \mu
_{2}^{^{\prime }}-mM\mu _{1}^{^{\prime }}\right) ^{2}}{\left( m+M\right) \mu
_{1}^{^{\prime }}-\mu _{2}^{^{\prime }}-mM},
\end{equation*}%
\begin{equation*}
\left( \mu _{4}^{^{\prime }}-\mu _{2}^{\prime ^{2}}\right) \left( \mu
_{2}^{^{\prime }}-\mu _{1}^{\prime ^{2}}\right) -\left( \mu _{3}^{^{\prime
}}-\mu _{1}^{^{\prime }}\mu _{2}^{\prime }\right) ^{2}\leq \frac{\left(
M-m\right) ^{6}}{432}
\end{equation*}%
and%
\begin{equation*}
\mu _{4}^{^{\prime }}-\mu _{2}^{\prime ^{2}}-\frac{\left( \mu _{3}^{^{\prime
}}-\mu _{1}^{^{\prime }}\mu _{2}^{\prime }\right) ^{2}}{\mu _{2}^{^{\prime
}}-\mu _{1}^{\prime ^{2}}}\leq \frac{\left( M-m\right) ^{4}}{64}.
\end{equation*}%
The inequalities (1.5) and (2.1) respectively give the upper bound for $\mu
_{2}$ and $\mu _{4}$ in terms of the range of the random variable, $M-m$. It
is interesting to note that the analogous upper bound for the third central
moment $\mu _{3}$ follows easily from the inequality (1.10).

\vskip0.1in\textbf{Theorem 2.3. }Let $X$ be a discrete or continuous random
variable taking values in\textbf{\ }$\left[ m,M\right] .\ $Then%
\begin{equation}
\left\vert \mu _{3}\right\vert \leq \frac{\left( M-m\right) ^{3}}{6\sqrt{3}}.
\tag{2.19}
\end{equation}%
\vskip0.1in\textbf{Proof. }From the inequality (1.10), we have%
\begin{equation}
\mu _{3}^{2}\leq \left( M-m\right) ^{2}\sigma ^{4}-4\sigma ^{6}.  \tag{2.20}
\end{equation}%
The inequality (2.19) follows from (2.20) and the fact that the function%
\begin{equation}
h\left( x\right) =\left( M-m\right) ^{2}x^{4}-4x^{6},  \tag{2.21}
\end{equation}%
achieves its maximum at $x=\frac{M-m}{\sqrt{6}}$ where $h\left( x\right)
\leq \frac{\left( M-m\right) ^{6}}{108}.\ \ \blacksquare $

\vskip0.1inEquality holds in (2.19) for $n=2;\ x_{1}=m\ $and $x_{2}=M$ with$%
\ p_{1}=\frac{1}{2}\pm \frac{1}{2\sqrt{3}}\ $and $p_{2}=\frac{1}{2}\mp \frac{%
1}{2\sqrt{3}}.$

It remains to prove an analogous of the Nagy inequality (1.6) for fourth
central moment. We show that a generalization of the Nagy inequality (1.6)
follows easily for the central moment $m_{2r}.$

\vskip0.1in\textbf{Theorem 2.4.} Let $m_{2r}$ be the central moment of $n$
real numbers $x_{i}$ such that $m\leq x_{i}\leq M,$ then%
\begin{equation}
m_{2r}\geq \frac{\left( M-m\right) ^{2r}}{2^{2r-1}n}+\left( \frac{n}{n-2}%
\right) ^{r-1}\left( m_{2}-\frac{\left( M-m\right) ^{2}}{2n}\right) ^{r}. 
\tag{2.22}
\end{equation}%
\vskip0.1in\textbf{Proof. }From (1.4), we have 
\begin{equation}
m_{2r}=\frac{\left( M-m_{1}^{^{\prime }}\right) ^{2r}+\left( m_{1}^{^{\prime
}}-m\right) ^{2r}}{n}+\frac{n-2}{n}\left( \frac{1}{n-2}\dsum%
\limits_{i=2}^{n-1}\left( x_{i}-m_{1}^{^{\prime }}\right) ^{2r}\right) . 
\tag{2.23}
\end{equation}%
It is evident that for $m$ positive real numbers $y_{1},y_{2,...}y_{m},$ 
\begin{equation}
\frac{1}{m}\dsum\limits_{i=1}^{m}y_{i}^{k}\geq \left( \frac{1}{m}%
\dsum\limits_{i=1}^{m}y_{i}\right) ^{k},\text{ }k=1,2,\ldots \text{ .} 
\tag{2.24}
\end{equation}%
Apply (2.24) to $n-2$ positive real numbers $\left( x_{i}-m_{1}^{^{\prime
}}\right) ^{2},i=2,\ldots ,n-1,$ we get%
\begin{equation}
\frac{1}{n-2}\dsum\limits_{i=2}^{n-1}\left( x_{i}-m_{1}^{^{\prime }}\right)
^{2r}\geq \left( \frac{1}{n-2}\dsum\limits_{i=2}^{n-1}\left(
x_{i}-m_{1}^{^{\prime }}\right) ^{2}\right) ^{r}.  \tag{2.25}
\end{equation}%
We also have%
\begin{equation}
\sum_{2}^{n-1}\left( x_{i}-m_{1}^{^{\prime }}\right) ^{2}=nm_{2}-\left(
m-m_{1}^{^{\prime }}\right) ^{2}-\left( M-m_{1}^{^{\prime }}\right) ^{2}. 
\tag{2.26}
\end{equation}%
Combine (2.23), (2.25) and (2.26), we have%
\begin{equation}
m_{2r}\geq \frac{\left( M-m_{1}^{^{\prime }}\right) ^{2r}+\left(
m_{1}^{^{\prime }}-m\right) ^{2r}}{n}+\frac{1}{n\left( n-2\right) ^{r-1}}%
\left( nm_{2}-\left( m-m_{1}^{^{\prime }}\right) ^{2}-\left(
M-m_{1}^{^{\prime }}\right) ^{2}\right) ^{r}.  \tag{2.27}
\end{equation}%
The right hand side expression (2.27) is minimum at $m_{1}^{^{\prime }}=%
\frac{m+M}{2}$, and so (2.22) follows from (2.27). $\blacksquare $

\vskip0.1inThe inequality (2.22) provides a generalization of the Nagy
inequality (1.6),%
\begin{equation}
m_{2r}\geq \frac{\left( M-m\right) ^{2r}}{2^{2r-1}n}.  \tag{2.28}
\end{equation}%
When $n=2$, the inequality (2.28) becomes equality. For $n=3$, equality
holds when $x_{1}=m,\ x_{2}=x_{3}=...=x_{n-1}=\frac{m+M}{2}$ and$\ x_{n}=M$.
Also, for $r=2$ and $n=3$, the inequalities (1.6) and (2.28) give equal
estimates.

\section{Bounds on the spread of a matrix}

\setcounter{equation}{0}Let $\mathbb{M}(n)$ be the space of all $n\times n$
complex matrices. A linear functional $\varphi :\mathbb{M}(n)\rightarrow 
\mathbb{C}
$ is said to be positive if $\varphi \left( A\right) $ is non-negative
whenever $A$ is positive semidefinite. It is unital if $\varphi \left(
I\right) =1.$ For more details, see \cite{4}. Let $A=\left( a_{ij}\right) $
be an element of $\mathbb{M}(n)$ with eigenvalues $\lambda _{i},$ $\
i=1,2,...,n.$ The spread of $A$ is defined as 
\begin{equation}
\text{spd}\left( A\right) =\max_{i,j}\left\vert \lambda _{i}-\lambda
_{j}\right\vert .  \tag{3.1}
\end{equation}%
It is shown in \cite[6]{5} that how positive unital linear maps can be used
to derive many inequalities for the spread. Enhancing this technique, we
derive here some more inequalities for the positive unital linear functional
and obtain bounds for the spread of a Hermitian matrix.

Beginning with Mirsky \cite{9} several authors have obtained bounds for the
spread of a matrix $A$ in terms of the functions of its entries. Mirsky \cite%
{10} proves that for every Hermitian matrix $A,$%
\begin{equation}
\text{spd}\left( A\right) ^{2}\geq \max_{i\neq j}\left( \left(
a_{ii}-a_{jj}\right) ^{2}+4\left\vert a_{ij}\right\vert ^{2}\right) . 
\tag{3.2}
\end{equation}%
Barnes and Hoffman \cite{1} prove the following sharper bound,%
\begin{equation}
\text{spd}\left( A\right) ^{2}\geq \max_{i,j}\left( \left(
a_{ii}-a_{jj}\right) ^{2}+2\underset{k\neq i}{\sum }\left\vert
a_{ik}\right\vert ^{2}+\underset{k\neq j}{2\sum }\left\vert
a_{jk}\right\vert ^{2}\right) \ .  \tag{3.3}
\end{equation}%
One more inequality of our present interest is, see \cite{5},%
\begin{equation}
\text{spd}\left( A\right) ^{2}\geq 4\max_{j}\underset{k\neq j}{\sum }%
\left\vert a_{jk}\right\vert ^{2}.  \tag{3.4}
\end{equation}%
The inequalities (3.3) and (3.4) are independent. Bhatia and Sharma \cite[6]%
{5} have shown that such inequalities follow easily from the inequalities
for positive linear maps. We pursue this topic further and obtain bounds for
the spread in the following theorems.

\vskip0.1in\textbf{Theorem 3.1. }Let $\varphi :\mathbb{M}(n)\rightarrow 
\mathbb{C}
$ be a positive unital linear functional and $A$ be any Hermitian element of 
$\mathbb{M}(n).$ Then%
\begin{equation}
\varphi \left( B^{4}\right) \leq \frac{\text{spd}\left( A\right) ^{4}}{12} 
\tag{3.5}
\end{equation}%
and%
\begin{equation}
\varphi \left( B^{2}\right) \varphi \left( B^{4}\right) -\varphi \left(
B^{2}\right) ^{3}-\varphi \left( B^{3}\right) ^{2}\leq \frac{\text{spd}%
\left( A\right) ^{6}}{432},  \tag{3.6}
\end{equation}%
where $B=A-\varphi \left( A\right) I.$

\vskip0.1in\textbf{Proof. } Let $\lambda _{i},$ $i=1,2,...,n$ be the
eigenvalues of $A$.$\ $By the spectral theorem,%
\begin{equation*}
B=\underset{i=1}{\overset{n}{\sum }}\left( \lambda _{i}-\varphi \left(
A\right) \right) P_{i},
\end{equation*}%
where $\lambda _{i}-\varphi \left( A\right) $ are the eigenvalues of $B$ and 
$P_{i}$ the corresponding projections with $\underset{i=1}{\overset{n}{\sum }%
}P_{i}=I$, see \cite{4}. Then, for $r=1,2,...,$ we have 
\begin{equation}
B^{r}=\underset{i=1}{\overset{n}{\sum }}\left( \lambda _{i}-\varphi \left(
A\right) \right) ^{r}P_{i}.  \tag{3.7}
\end{equation}%
Apply $\varphi $ to both sides of (3.7), we get%
\begin{equation}
\varphi \left( B^{r}\right) =\underset{i=1}{\overset{n}{\sum }}\left(
\lambda _{i}-\varphi \left( A\right) \right) ^{r}\varphi \left( P_{i}\right)
.  \tag{3.8}
\end{equation}%
Since $\lambda _{i}-\varphi \left( A\right) $ are real numbers and $\varphi
\left( P_{i}\right) $ are non-negative real numbers such that $\underset{i=1}%
{\overset{n}{\sum }}\varphi \left( P_{i}\right) =1$, the inequalities (3.5)
and (3.6) follow respectively from (2.1) and (2.9). \ \ $\blacksquare $

Note that an equivalent form of (3.6) says that the determinant%
\begin{equation*}
\left\vert 
\begin{array}{ccc}
1 & \varphi \left( A\right) & \varphi \left( A^{2}\right) \\ 
\varphi \left( A\right) & \varphi \left( A^{2}\right) & \varphi \left(
A^{3}\right) \\ 
\varphi \left( A^{2}\right) & \varphi \left( A^{3}\right) & \varphi \left(
A^{4}\right)%
\end{array}%
\right\vert \leq \frac{\text{spd}\left( A\right) ^{6}}{432}.
\end{equation*}%
In this connection we prove one more inequality in the following theorem.

\vskip0.1in\textbf{Theorem 3.2.} Let $\varphi :\mathbb{M}(n)\rightarrow 
\mathbb{C}
$ be a positive unital linear functional. Then for $0\leq A\leq MI$, we have%
\begin{equation}
0\leq \left\vert 
\begin{array}{cc}
\varphi \left( A\right) & \varphi \left( A^{2}\right) \\ 
\varphi \left( A^{2}\right) & \varphi \left( A^{3}\right)%
\end{array}%
\right\vert \leq \frac{M^{4}}{27}.  \tag{3.9}
\end{equation}

\vskip0.1in\textbf{Proof.} On using arguments similar to those used in the
proof of the above theorem, it follows from the second inequality (1.9) that%
\begin{equation}
\varphi \left( A^{3}\right) \leq M\varphi \left( A^{2}\right) -\frac{\left(
M\varphi \left( A\right) -\varphi \left( A^{2}\right) \right) ^{2}}{%
M-\varphi \left( A\right) }.  \tag{3.10}
\end{equation}%
Since $\varphi \left( A\right) \geq 0$, the inequality (3.10) implies that
for $A<MI$, 
\begin{equation}
\varphi \left( A^{3}\right) \varphi \left( A\right) -\varphi \left(
A^{2}\right) ^{2}\leq M\varphi \left( A^{2}\right) \varphi \left( A\right) -%
\frac{\left( M\varphi \left( A\right) -\varphi \left( A^{2}\right) \right)
^{2}\varphi \left( A\right) }{M-\varphi \left( A\right) }-\varphi \left(
A^{2}\right) ^{2}.  \tag{3.11}
\end{equation}%
The inequality (3.9) follows from (3.11) and the fact that the function%
\begin{equation*}
h\left( x,y\right) \leq bxy-\frac{\left( bx-y\right) ^{2}}{b-x}x-y^{2},
\end{equation*}%
achieves its maximum at $x=\frac{2}{3}b$ \ and\ $y=\frac{x\left( b+x\right) 
}{2}$, where $h\left( x,y\right) \leq \frac{b^{4}}{27}$. $\ \ \blacksquare $

We now consider an upper bounds for the spread of a matrix. Mirsky \cite{9}
proves that for any $n\times n$ matrix $A$,%
\begin{equation*}
\text{spd}\left( A\right) ^{2}\leq 2\text{tr}A^{\ast }A-\frac{2}{n}%
\left\vert \text{tr}A\right\vert ^{2}.
\end{equation*}%
See also \cite[19]{2}. We prove an extension of this inequality in the
following theorem.

\vskip0.1in\textbf{Theorem 3.3. } Let $A$ be $n\times n$\ matrix, then%
\begin{equation}
\text{spd}\left( A\right) ^{2r}\leq 2^{2r-1}\text{tr}\left( B^{r}\left(
B^{\ast }\right) ^{r}\right) ,  \tag{3.12}
\end{equation}%
where $B=A-\frac{\text{tr}A}{n}I\ $and $r=1,2,...$

\vskip0.1in\textbf{Proof. }Let $\lambda _{i}$ be the eigenvalues of $A,\
i=1,2,...,n$. Then%
\begin{equation}
\frac{1}{n}\text{tr}\left( B^{r}\left( B^{\ast }\right) ^{r}\right) \geq 
\frac{1}{n}\sum_{i=1}^{n}\left\vert \lambda _{i}-\frac{\text{tr}A}{n}%
\right\vert ^{2r}.  \tag{3.13}
\end{equation}%
From (3.13), we see that the inequality%
\begin{equation}
\frac{1}{n}\text{tr}\left( B^{r}\left( B^{\ast }\right) ^{r}\right) \geq 
\frac{1}{n}\left( \left\vert \lambda _{j}-\frac{\text{tr}A}{n}\right\vert
^{2r}+\left\vert \lambda _{k}-\frac{\text{tr}A}{n}\right\vert ^{2r}\right) ,
\tag{3.14}
\end{equation}%
holds for any $j,k=1,2,...n$ with $j\neq k.$\ Also, for two positive real
numbers $x_{1}$and $x_{2}$, $2^{r-1}\left( x_{1}^{r}+x_{2}^{r}\right) \geq
\left( x_{1}+x_{2}\right) ^{r},$ therefore%
\begin{equation}
\left\vert \lambda _{j}-\frac{\text{tr}A}{n}\right\vert ^{2r}+\left\vert
\lambda _{k}-\frac{\text{tr}A}{n}\right\vert ^{2r}\geq \frac{1}{2^{2r-1}}%
\left( \left\vert \lambda _{j}-\frac{\text{tr}A}{n}\right\vert +\left\vert
\lambda _{k}-\frac{\text{tr}A}{n}\right\vert \right) ^{2r}.  \tag{3.15}
\end{equation}%
Using triangular inequality, we have%
\begin{equation}
\left\vert \lambda _{j}-\frac{\text{tr}A}{n}\right\vert +\left\vert \lambda
_{k}-\frac{\text{tr}A}{n}\right\vert \geq \left\vert \lambda _{j}-\lambda
_{k}\right\vert .  \tag{3.16}
\end{equation}%
Combine (3.14)-(3.16); we easily get (3.12). \ \ $\blacksquare $

Several inequalities for the spread can be obtained from (3.5) and (3.6).
For example, if we choose $\varphi \left( A\right) =\frac{\text{tr}A}{n},$
we have%
\begin{equation}
\text{spd}\left( A\right) ^{4}\geq \frac{12}{n}\text{tr}B^{4}  \tag{3.17}
\end{equation}%
and%
\begin{equation}
\text{spd}\left( A\right) ^{6}\geq \frac{432}{n^{3}}\left( n\text{tr}B^{2}%
\text{tr}B^{4}-\left( \text{tr}B^{2}\right) ^{3}-n\left( \text{tr}%
B^{3}\right) ^{2}\right) .  \tag{3.18}
\end{equation}%
Note that (1.8) yields first inequality (1.7) for $\varphi \left( A\right) =%
\frac{\text{tr}A}{n}$. Also, (1.8) gives (3.3) and (3.4) respectively for $%
\varphi \left( A\right) =\frac{a_{ii}+a_{jj}}{2}$ and $\varphi \left(
A\right) =a_{ii}$. The corresponding estimates for the spread from (3.5) and
(3.6) can be calculated numerically, see Example 1, below.

We give examples and compare the bound (1.7) in terms of the traces with our
corresponding bounds (3.12), (3.17) and (3.18). Likewise, we compare
(3.3)-(3.4) with (3.5)-(3.6), respectively.

\vskip0.1in\textbf{Example 1. }Let

\begin{equation*}
A=\left[ 
\begin{array}{ccc}
3 & 2 & 1 \\ 
2 & 0 & 2 \\ 
1 & 2 & 3%
\end{array}%
\right] .
\end{equation*}%
Then from the bound (1.7), spd$\left( A\right) \geq 5.6569$ while from our
bounds (3.17) and (3.18) we respectively have spd$\left( A\right) \geq
\allowbreak 5.8259$ and spd$\left( A\right) \geq 6.9282$. Here $n=3$, the
inequalities (1.7) and (3.12) therefore give equal estimates spd$\left(
A\right) \leq \allowbreak 6.9282$, $r=2$. Further, from (3.3) spd$\left(
A\right) \geq \allowbreak 5.9161$ while from our bounds (3.5) and (3.6) for $%
\varphi \left( A\right) =\frac{a_{ii}+a_{jj}}{2}$ give spd$\left( A\right)
\geq \allowbreak 6.\,\allowbreak 0181$ and spd$\left( A\right) \geq
\allowbreak 6.8252$. Likewise, from (3.4), spd$\left( A\right) \geq 5.6569$
and from (3.5) and (3.6) we respectively have spd$\left( A\right) \geq
\allowbreak 6.\,\allowbreak 9282$ and spd$\left( A\right) \geq \allowbreak
6.2947$, $\varphi \left( A\right) =a_{ii}$. So, our bounds give better
estimates.

\textbf{Example 2. }Let%
\begin{equation*}
A_{1}=\left[ 
\begin{array}{cccc}
6 & 3 & 4 & 2 \\ 
3 & 1 & 0 & 3 \\ 
4 & 0 & 2 & 1 \\ 
2 & 3 & 1 & 2%
\end{array}%
\right] \ \ \text{and\ \ }A_{2}=\left[ 
\begin{array}{cccc}
6 & 0 & 4 & 2 \\ 
3 & 1 & 0 & 3 \\ 
4 & 0 & 2 & 1 \\ 
2 & 3 & 1 & 2%
\end{array}%
\right] .
\end{equation*}%
For the Hermitian matrix $A_{1}$, (1.7) gives spd$\left( A_{1}\right) \leq $ 
$\allowbreak 13.\,\allowbreak 620$ while from our bound (3.12) spd$\left(
A_{1}\right) \leq $ $\allowbreak 13.559,\ r=2$. Likewise, for arbitrary
matrix $A_{2}$ the Mirsky bound (3.12) with $r=1$ gives spd$\left(
A_{2}\right) \leq \allowbreak 12.\,\allowbreak 227$ while from our bound
(3.12), spd$\left( A_{2}\right) \leq \allowbreak 11.934,\ r=2$.

\section{Bounds for the span of a polynomial}

\setcounter{equation}{0}In the theory of polynomial equations, the study of
polynomials with real roots is of special interest, see \cite[14]{8}. The
span of a polynomial is the length $b-a$ of the smallest interval $\left[ a,b%
\right] $ containing all the zeros of polynomial. It is also of interest to
find bounds on the roots and span of a polynomial in terms of its
coefficients; see \cite[14, 15]{7}. We obtain here some bounds for the span
of polynomial.

It is sufficient to consider the polynomial equation in which the
coefficient of $x^{n-1}$ is zero,%
\begin{equation}
f\left( x\right) =x^{n}+a_{2}x^{n-2}+a_{3}x^{n-3}+\ldots +a_{n-1}x+a_{n}=0%
\text{.}  \tag{4.1}
\end{equation}%
Let $x_{1},x_{2},...,x_{n}$ be the roots of (4.1). On using the well known
Newton's identity 
\begin{equation*}
\alpha _{k}+a_{1}\alpha _{k-1}+a_{2}\alpha _{k-2}+\ldots +a_{k-1}\alpha
_{1}+ka_{k}=0,\text{ }
\end{equation*}%
where $\alpha _{k}=\sum\limits_{i=1}^{n}x_{i}^{k}$ and $k=1,2,\ldots n,$ we
have 
\begin{equation}
m_{1}=\frac{1}{n}\sum_{i=1}^{n}x_{i}=0,\text{ }m_{2}=\frac{1}{n}%
\sum_{i=1}^{n}x_{i}^{2}=-\frac{2}{n}a_{2},\text{ }m_{3}=\frac{1}{n}%
\sum_{i=1}^{n}x_{i}^{3}=-\frac{3}{n}a_{3}  \tag{4.2}
\end{equation}%
and%
\begin{equation}
m_{4}=\frac{1}{n}\sum_{i=1}^{n}x_{i}^{4}=\frac{2}{n}\left(
a_{2}^{2}-2a_{4}\right) .  \tag{4.3}
\end{equation}%
The span of polynomial (4.1) is spn$\left( f\right) =\underset{i,j}{\max }%
\left\vert x_{i}-x_{j}\right\vert $. Then, from (1.6), we get the Nagy's
inequality, see \cite[14]{11}, 
\begin{equation}
\text{spn}\left( f\right) \leq 2\sqrt{-a_{2}}.  \tag{4.4}
\end{equation}%
Likewise, from (1.5) we have%
\begin{equation}
\text{spn}\left( f\right) \geq 2\sqrt{\frac{-2a_{2}}{n}},  \tag{4.5}
\end{equation}%
see \cite{14}.

In a similar spirit, we obtain some further estimates for spn$\left(
f\right) $ in the following theorems.

\vskip0.1in\textbf{Theorem 4.1. }If the roots of the polynomial\textbf{\ }%
(4.1) are all real, then for $n\geq 5,$ we have 
\begin{equation}
\left( \frac{24}{n}\left( a_{2}^{2}-2a_{4}\right) \right) ^{\frac{1}{4}}\leq 
\text{spn}\left( f\right) \leq 2\left( a_{2}^{2}-2a_{4}\right) ^{\frac{1}{4}%
}.  \tag{4.6}
\end{equation}%
\vskip0.1in\textbf{Proof. }Let $x_{i}$ be the roots of polynomial (4.1) such
that $x_{1}\leq x_{i}\leq x_{n},\ i=1,2,\ldots n$. Then from the inequality
(2.1), we have%
\begin{equation}
\left( x_{n}-x_{1}\right) ^{4}\geq 12m_{4}.  \tag{4.7}
\end{equation}%
Combine (4.3) and (4.7), we immediately get the first inequality (4.6), $%
x_{n}-x_{1}=$spn$\left( f\right) $. Similarly, the inequality (2.28) gives $%
x_{n}-x_{1}\leq \left( 8nm_{4}\right) ^{\frac{1}{4}},\ r=2$. This implies
the second inequality (4.6). $\ \ \blacksquare $

\vskip0.1in\textbf{Theorem 4.2. }Under the conditions of Theorem 4.1, we have%
\begin{equation}
\text{spn}\left( f\right) \geq \left( \frac{432}{n^{3}}\left( 4\left(
2-n\right) a_{2}^{3}-9na_{3}^{2}+8na_{2}a_{4}\right) \right) ^{\frac{1}{6}}.
\tag{4.8}
\end{equation}

\vskip0.1in\textbf{Proof. }As in the proof of above theorem, it follows from
the inequality (2.9) that%
\begin{equation}
\left( x_{n}-x_{1}\right) ^{6}\geq 432\left(
m_{2}m_{4}-m_{2}^{3}-m_{3}^{2}\right) .  \tag{4.9}
\end{equation}%
Combine (4.2), (4.3) and (4.9), we immediately get (4.8).\ $\ \ \blacksquare 
$

\vskip0.1in\textbf{Example. }Let%
\begin{equation}
f\left( x\right) =x^{5}+80x^{4}+1500x^{3}+5000x^{2}+3750x+\frac{1}{5}=0. 
\tag{4.10}
\end{equation}%
The roots $x_{i}$ of (4.10) are real, $i=1,2,...,5;\ $see \cite{18}. Let $%
y_{i}=x_{i}-16$ be the roots of the diminished equation 
\begin{equation*}
f\left( y\right) =y^{5}-1060y^{3}+14920y^{2}+12710y-\frac{3648479}{5}=0.
\end{equation*}%
The Nagy inequality (4.4) gives $s\left( f\right) $ $\leq 65.116$ while from
(4.6) $s\left( f\right) \leq 64.744$. Also, the Popoviciu inequality (4.5)
gives $s\left( f\right) \geq 41.183$ while from our bounds (4.6) and (4.8)
we respectively have $s\left( f\right) \geq 47.916\ $and $s\left( f\right)
\geq 48.435.$

\vskip0.2in\textbf{Acknowledgements}. The authors are grateful to Prof.
Rajendra Bhatia for the useful discussions and suggestions. The first two
authors thank I.S.I. Delhi for a visit in January 2014 when this work had
begun.

\end{document}